\documentclass{article}

\usepackage[english]{babel}
\usepackage[utf8]{inputenc}
\usepackage[T1]{fontenc}

\usepackage[a4paper,top=3cm,bottom=2cm,left=3cm,right=3cm,marginparwidth=1.75cm]{geometry}

    \usepackage{amsthm}
    \usepackage{amssymb}
    \usepackage{amsmath}
    \usepackage{mathtools}
    \usepackage{epsfig}
    \usepackage{caption}
    \usepackage[colorinlistoftodos]{todonotes}
    \usepackage[colorlinks=true, allcolors=blue]{hyperref}
    \usepackage{xcolor}
    \usepackage{pinlabel}
    \usepackage{geometry}
    \usepackage{subcaption}
    \usepackage[numbers]{natbib} 
    \usepackage{hyperref}
    \usepackage{authblk}
    \usepackage{longtable}
    \usepackage[shortlabels]{enumitem}

\theoremstyle{definition}

\theoremstyle{plain}

\theoremstyle{remark}

\newtheorem{example}{Example}[section]

\begin{document}

\title{A Markov theorem for plat closure of surface braids in Dunwoody and periodic Takahashi manifolds}

\author[1]{Alessia Cattabriga}
\author[1]{Paolo Cavicchioli}

\affil[1]{Alma Mater Studiorum, Bologna, Italy}

\date{February 2022}

\maketitle

\begin{abstract}
In this article we deal with the problem of finding equivalence moves for links in Dunwoody and periodic Takahashi manifolds. We represent these manifolds using Heegaard splitting and we represent the embedded links as plat closure of elements in the braid group of the corresponding Heegaard surfaces. More precisely, starting from an open Heegaard diagram for such manifolds, we determine the plat slide equivalence moves algorithmically and compute them explicitly in some cases.
\end{abstract}



\section{Introduction} 

The representation of links as plat closures of braids with an even number of strands goes back to the works of Hilden \cite{hilden1975generators} and Birman \cite{birman1976stable}, who proved that any link in $\mathbb{R}^3$ may be represented as the plat closure of a braid. Since then, much work has been done using this representation in the direction of studying the equivalence problem of links, for example by defining and analyzing link invariants as the  Jones polynomial (see \cite{bigelow2002does, birman1988jones}).
Using Heegaard splittings, in \cite{doll1993generalization} Doll introduced the notion of $(g,b)$-decomposition or generalized bridge decomposition for links in a closed, connected and orientable (from now on c.c.o.) \(3\)-manifold, opening the way to the study of links in 3-manifolds via surface braid groups and their plat closures \cite{bellingeri2012hilden, cattabriga2004mulazzani, tawn2008}. 
The equivalence of links in 3-manifolds, under isotopy, has recently been described in \cite{cattabriga2018markov}, where the authors find a finite set of moves that connect braids in \(B_{g, 2n}\), the braid group on \(2n\) strands of a surface of genus $g$, having isotopic plat closures. 
In their result some  moves are explicitly described as elements in \(B_{g, 2n}\), while others, called \textit{plat slide moves}, depend on a Heegaard surface for the manifold and are explicitly described only for the case of Heegaard genus one, that is for lens spaces and \(S^2 \times S^1\). Moreover, in \cite{cavicchioli2021algorithmic} an algorithm to determine  plat slide moves  for manifolds with Heegaard genus two is presented.  In this paper we deal with the same problem, but in the case of two infinite families of manifolds: Dunwoody and (periodic) Takahashi manifolds.

Dunwoody manifolds are c.c.o. 3-manifolds introduced in \cite{dunwoody1995cyclic} by means of a graph with a cyclic symmetry determining an open Heegaard diagram. Some interesting results connect these manifolds and cyclic branched covering of a \((1, 1)\)-knot (i.e., knots admitting a \((1,1)\) decomposition): in \cite{cattabriga2004all} and \cite{grasselli2001genus} it has been proven that the class of Dunwoody manifolds coincides with the class of strongly-cyclic branched coverings of \((1, 1)\)-knots. In \cite{cattabriga2010mulazzanivesnin} the family of Dunwoody manifold has been generalized to include also manifolds with non-empty boundary, while, more recently,  other graph with cyclic symmetry generalizing those introduced in  \cite{dunwoody1995cyclic} have been studied (see \cite{howie2020williams}).

Takahashi manifolds are c.c.o. 3-manifolds introduced in \cite{takahashi1989presentations} by Dehn surgery with rational coefficients  along a specific $2n$-component link in $S^3$ (see Figure \ref{fig:takahashi_dehn_basic}). 
Several important classes of 3-manifolds, such as (fractional) Fibonacci manifolds \cite{helling1998geometric, vesnin1998fractional} and Sieradski manifolds \cite{cavicchioli1998geometric}, represent notable examples of periodic Takahashi manifolds. These have been intensively studied in many papers as \cite{kim2000class, mulazzani2001periodic, ruini1998structure, vesnin1998fractional}. Among Takahashi manifolds the periodic ones are those whose surgery coefficients  present an order $n$  cyclic symmetry: in \cite{mulazzani2001periodic} the author proves that these manifolds are cyclic coverings of the connected sum of two lens spaces or $S^3$ branched over a $(2,1)$-knot.

In the paper, we determine algorithmically plat slide moves for links in these manifolds: in order to do so, we use the fact that both families admit a symmetric open Heegaard diagram depending of a finite number of integer parameters. Our approach could be easily generalized  to other families of manifolds with the same property. 


The paper is organized as follows. 
In Section \ref{Preliminaries} we recall the definition of Heegaard diagram, plat closure of links in 3-manifolds and the equivalence moves introduced in \cite{cattabriga2018markov}. In Section \ref{Dunwoody} (resp. \ref{Takahashi}), after giving the definition of the family of Dunwoody (resp. Takahashi) manifolds, we compute the equivalence moves for the considered family of 3-manifolds, together with some notable examples. 

In the following, manifolds are always assumed to be closed, connected and orientable and are considered up to homemorphisms. Links inside them, that is, closed 1-dimensional smooth submanifolds, will be considered up to isotopy.

\section{Preliminaries} \label{Preliminaries}

In this section we briefly recall the notion of Heegaard splitting and Heegaard diagram for 3-manifolds as well as the definition of plat closure of links in 3-manifolds.

\subsection{Heegaard splittings}\label{prel_heeg}

\begin{figure}[h!]
    \centering
    \includegraphics[width = .5\textwidth]{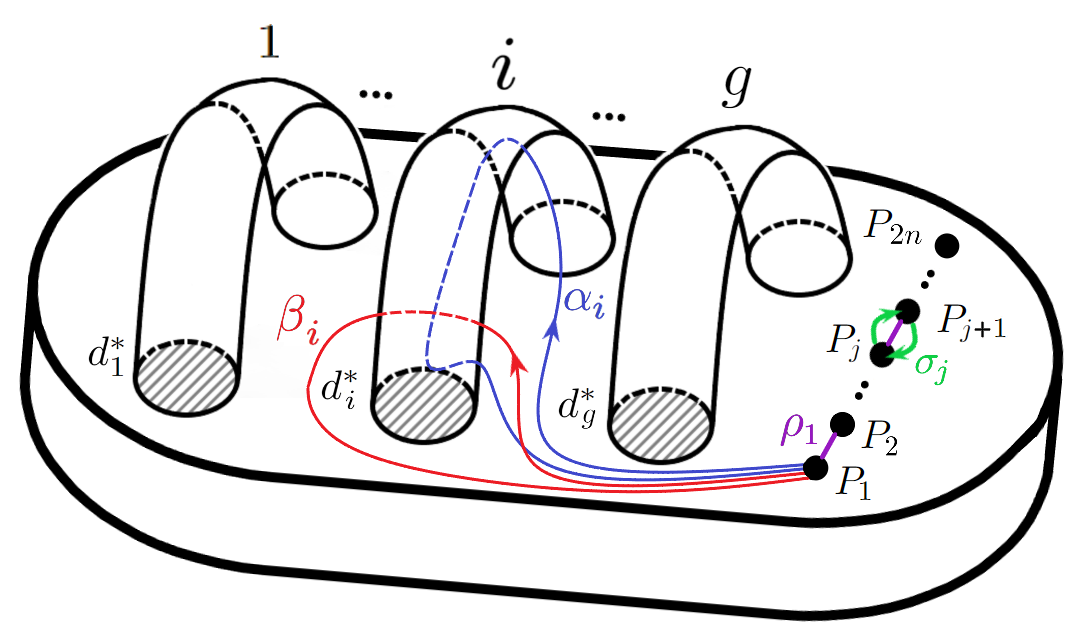}
    \caption{Standard genus $g$ handlebody and generators of the surface braid group.}
    \label{fig:generator_dunwoody}
\end{figure}

Let \(M\) be a 3-manifold, a \emph{Heegaard spitting} for $M$ is the data of a triple \((H^*,H,\phi)\), where $H^*$ and $H$ are two copies of a  genus $g$ oriented handlebody (see Figure \ref{fig:generator_dunwoody}) standardly embedded in $\mathbb R^3$ and  \(\phi:\partial H \rightarrow \partial H^*\) is an orientation reversing homeomorphism such that $M=H \cup_{\phi} H^*$. The surface \(\partial H\cup_{\phi}\partial H^*\subseteq M\) has genus $g$ and  is called \emph{Heegaard surface} for $M$. 

Each \(3\)-manifold admits Heegaard splittings \cite{heegaard1898forstudier}: we call \emph{Heegaard genus} of a 3-manifold \(M\) the minimal genus of a Heegaard surface for \(M\). \\ 
For instance, the \(3\)-sphere \(S^3\) is the only \(3\)-manifold with Heegaard genus \(0\), while the manifolds with Heegaard genus \(1\) are lens spaces (i.e., cyclic quotients of \(S^3\)) and \(S^2\times S^1\). While \(S^3\), as well as lens spaces, have, up to isotopy, only one Heegaard surface of minimal genus (and those of higher genera are stabilizations of that of minimal genus), in general, a manifold may admit non isotopic Heegaard surfaces of the same genus (see \cite{moriah1988heegaard} for the case of Seifert manifolds).\\
Each Heegaard splitting  \((H^*,H,\phi)\) can be represented by means of a (closed) \emph{Heegaard diagram} that is a triple $(\Sigma,\mathbf{d}^*,\mathbf{d})$, with $\Sigma=\partial H^*$ and $\mathbf{d}^*= \{d_1^*, \dots, d_g^*\}$ (resp. $\mathbf{d}= \{d_1, \dots, d_g\}$) a set of meridians\footnote{We recall that a \emph{set of meridians} for a handlebody of genus \(g\) is a collection of closed curves \(\mathbf{d}^* = \{d^*_1, \dots, d^*_g\}\) on its boundary, such that the curves bound properly embedded pairwise disjoint disks \(D_1, \dots, D_g\) (dotted in  Figure \ref{fig:generator_dunwoody}) and cutting the handlebody along these disks yields a 3-ball.} for $H^*$ (resp. the image in $\Sigma$, via $\phi$, of a set of meridians for $H$). 
Moreover, one of the two system of meridians, say $\mathbf{d}^*$ could be always choose in a standard way, as in Figure \ref{fig:generator_dunwoody}. 
If we cut \(\Sigma\) along \(\mathbf{d}^*\), we obtain a sphere with \(2g\) holes, say \(D_1, \dots, D_{2g}\), distinct and paired, each pair corresponding to a certain meridian \(d_i^*\). Moreover the curves of  \(\mathbf{d}\) will be naturally cut into arcs joining the holes in various ways, giving a graph on the sphere with \(2g\) holes called \emph{open Heegaard diagram} of  \((H^*,H,\phi)\) and still denoted by \((\Sigma, \mathbf{d}^*, \mathbf{d})\). \\
On the contrary, if we have an open Heegaard diagram, in order to obtain the closed one we just need to identify the \(D_i\)'s. It can be done by attaching 1-handles connecting the paired disks and arcs along the handles, connecting the paired vertices of the graph.\\ 
The result will be a set of closed curves which will represent the system \(\mathbf{d}\). 

\subsection{Plat closure of links in 3-manifolds}\label{prel_links}

Let $M$ be a 3-manifold, fix a Heegaard diagram $(\Sigma,\mathbf{d}^*,\mathbf{d})$ for $M$ so that $\Sigma$ is the boundary of a standard handlebody and $\mathbf{d}^*$ are chosen in a standard way (see Figure \ref{fig:generator_dunwoody}); denote with $g$ the genus of $\Sigma$. Consider $B_{g,2n}=\pi_1(C_{2n}(\Sigma), *)$ the braid group on $2n$ strands of $\Sigma$, that is, the fundamental group of the configuration space of $2n$ (unordered)  points $P_1,\ldots, P_{2n}$ in $\Sigma$. A presentation  of   $B_{g,2n}$  is given in \cite{bellingeri2004presentations}. 
Referring to Figure \ref{fig:generator_dunwoody} the generators are   $\sigma_1,\ldots,\sigma_{2n-1}$, the standard braid ones, and $\alpha_1,\ldots,\alpha_g,\beta_1,\ldots, \beta_g$, where $\alpha_i$ (resp. $\beta_i$) is the braid whose strands are all trivial except the first one which goes once along the $i$-th longitude (resp. $i$-th meridian) of $\Sigma$.  Moreover, we fix a set of \(n\) disjoint arcs  \(\rho_1,\ldots,\rho_n\) embedded into \(\Sigma\), such that   \(\partial \rho_i=\{P_{2i-1},P_{2i}\}\), for \(i =1,\ldots,n\), as in Figure \ref{fig:generator_dunwoody}. Now we can associate to each element \(\gamma\in B_{g, 2n}\) a link in $M$, called  the  \emph{plat closure} of $\gamma$, and denoted with \(\hat{\gamma}\), obtained ``closing'' a geometric representative of \(\gamma\) in $\Sigma\times [0,1]$ by connecting   \(P_{2i-1}\times\{0\}\) with  \(P_{2i}\times\{0\}\) through   \(\rho_i\times\{0\}\) and \(P_{2i-1}\times\{1\}\) with  \(P_{2i}\times\{1\}\)  through  \(\rho_i\times\{1\}\), for \(i=1,\ldots,n\).

For each link \(L\) in \(M\), there exists a braid \(\gamma \in \cup_{n \in \mathbb{N}} B_{g, 2n}\) such that \(\hat{\gamma}=L\). 
Moreover, as proved  in  \cite[Theorem 3]{cattabriga2018markov}, two braids \(\gamma_1, \gamma_2 \in \bigcup_{n\in\mathbb N} B_{g,2n}\) have isotopic plat closures if and only if \(\gamma_1\) and \(\gamma_2\) differ by a finite sequence of the following moves
\begin{align*}
  (M1)  & \qquad \sigma_1\gamma \longleftrightarrow  \gamma  \longleftrightarrow \gamma \sigma_1\\
  (M2)  & \qquad \sigma_{2i}\sigma_{2i+1}\sigma_{2i-1}\sigma_{2i} \gamma \longleftrightarrow \gamma \longleftrightarrow \gamma \sigma_{2i}\sigma_{2i+1}\sigma_{2i-1}\sigma_{2i} \\
  (M3)  & \qquad \sigma_2\sigma_1^2\sigma_2\gamma \longleftrightarrow  \gamma  \longleftrightarrow \gamma \sigma_2\sigma_1^2\sigma_2\\
   (M4)  & \qquad \alpha_j\sigma_1^{-1}\alpha_j\sigma_1^{-1}\gamma \longleftrightarrow \gamma\longleftrightarrow \gamma \alpha_j\sigma_1^{-1}\alpha_j\sigma_1^{-1}\quad \textup{for } j=1,\ldots, g \\
  (M5) &\qquad \beta_j\sigma_1^{-1}\beta_j\sigma_1^{-1}\gamma \longleftrightarrow \gamma\longleftrightarrow \gamma \beta_j\sigma_1^{-1}\beta_j\sigma_1^{-1}\quad \textup{for } j=1,\ldots, g \\
  (M6) &\qquad \gamma \longleftrightarrow T_k(\gamma)\sigma_{2k}\\
  \ &\ \ \quad \textup{ where } T_k: B_{g,2n}\to B_{g, 2n+2} \textup{ is  defined by } 
   T_k(\alpha_i)=\alpha_i, \ T_k(\beta_i)=\beta_i \textup{ and }\\
\ & \ \ \quad T_k(\sigma_i)=\left\{\begin{array}{l}\sigma_i \qquad \qquad \qquad \qquad \qquad \qquad\quad\textup{if } i<2k\\
\sigma_{2k}\sigma_{2k+1}\sigma_{2k+2}\sigma_{2k+1}^{-1}\sigma_{2k}^{-1}\ \qquad \quad \textup{if } i=2k\\
\sigma_{i+2} \qquad \qquad \qquad \qquad \qquad \qquad \textup{if } i>2k\\
\end{array}\right.\\
(psl_i^*) &\qquad \gamma \longleftrightarrow  \gamma \beta_i \quad  \textup{for } i = 1, \dots, g\\
(psl_i) &\qquad \gamma \longleftrightarrow \bar{d}_i \gamma, \quad  \textup{for } i = 1, \dots, g,
\end{align*}
where \(\bar{d}_i\in B_{g,2}\) is a braid representative of \(d_i\) with $2$ strands.\footnote{Since $d_i$ bounds a disk in $M$ it is always possible such a representation.}
So in order to determine completely the moves, it is necessary to describe the \(\bar{d}_i\) as a word in terms of the braid generators. 
From now on we will call \textit{psl-moves} the moves \(psl_i\). 

\subsection{From Heegaard diagrams to \textit{psl}-moves}

In \cite{cattabriga2018markov} an explicit formula for \textit{psl}-moves in 3-manifold with Heegaard genus one is computed, while in  \cite{cavicchioli2021algorithmic}, an algorithm to compute the \textit{psl}-moves for 3-manifolds of Heegaard genus two is described, starting from an open Heegaard diagram of the manifold. However, the same idea works in higher genus, so we briefly recall it. 

Let \((\Sigma, \mathbf{d}^*, \mathbf{d})\) be an open Heegaard diagram for $M$. In the corresponding closed one on $\Sigma$, we color in blue the arcs along the 1-handles arising from the identifications of the paired disks, and in red the other arcs. Note that, by construction, moving along each $d_i$, any red arc is followed by a blue arc and vice-versa.

In order to represent \(d_i\)'s curves in terms of generators of the braid group we can proceed as follows: (1) orient each curve and divide it into elementary pieces of couples of subsequent red-blue arcs, (2) connect the two ends of each piece with a fixed base point by fixed arcs, so to obtain an elementary closed curve, (3) catalogue each elementary closed curve arising in this way up to homotopy and (4) represent it in terms of generators of the surface braid group. This procedure gives rise to a \textit{dictionary} that can be used to implement a computer algorithm that takes as input the graph corresponding to the open Heegaard diagram and gives as output the words $\bar{d}_i$'s corresponding to the \textit{psl}-moves. It just works by visiting the cycles of the graph and replacing each elementary arc with the corresponding word in the dictionary.

\section{Dunwoody manifolds} \label{Dunwoody}
\begin{figure}[h!]
    \centering
    \includegraphics[width = .7\textwidth]{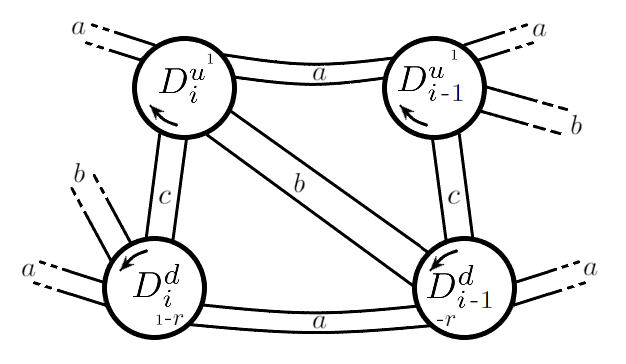}
    \caption{The graph \(\Gamma(a,b,c,n)\) and the labelling of its vertices. An arc with label \(k\) identifies \(k\) parallel arcs.}
    \label{fig:dunwoody_label}
\end{figure}

Let \(a, b, c, n \in \mathbb{N}\), with \(n>0\) and \(a+b+c>0\). Let \(\Gamma = \Gamma(a,b,c,n)\) be the regular, trivalent, planar graph depicted in  Figure \ref{fig:dunwoody_label}. The graph contains \(n\) top and \(n\) bottom circles, each having \(d = 2a+b+c\) vertices; we denote them with \(D_1^u, \dots, D_n^u\) and \(D_1^d, \dots, D_n^d\) respectively (all the indexes are considered mod \(n\)).  For each \(i = 1, \dots, n\), there are \(a\) parallel arcs, called \textit{upper} (resp. \textit{lower}) arcs, connecting \(D_{i-1}^u\) (resp. \(D_{i-1}^d\)) with \(D_{i}^u\) (resp. \(D_{i}^d\)), \(b\) parallel arcs, called \textit{diagonal} arcs, connecting \(D_i^u\) to \(D_{i-1}^d\) \(c\) parallel arcs, called \textit{vertical} arcs, connecting \(D_i^u\) to \(D_i^d\). Clearly, $\Gamma$ has a cyclic symmetry of order $n$ sending  \(D_i^u\) into \(D_{i-1}^u\) and \(D_i^d\) into \(D_{i-1}^d\), for $i=1,\ldots,n$.

The one-point compactification of the plane brings \(\Gamma\) onto \(S^2\).   Now, let \(r, s\) two given integers: we give a clockwise (resp. counterclockwise) orientation to each \(D_i^u\) (resp. \(D_i^d\)) and we enumerate the vertices as in Figure \ref{fig:dunwoody_label} (all the numbers are to be considered mod \(d\)), for \(i = 1, \dots, n\).  If we cut of from \(S^2\) the region of the \(2n\) disks bounded by each \(D_i^u\) and \(D_i^d\) and not containing any arc of \(\Gamma\), and we glue  \(D_i^u\) with \(D_{i+s}^d\), for $i=1,\ldots,n$, so that vertices having the same labelling correspond,  we obtain a regular graph of degree four in an orientable surface \(\Sigma_n\) of genus \(n\). Of course, by construction, we can always consider \(r\) mod \(d\) and \(s\) mod \(n\). 

Through the identification, the \(nd\) arcs of \(\Gamma\) connect with each other via their endpoints creating \(m\) closed curves \(e_1, \dots, e_m\), where \(e_i\) is the curve passing from the vertex labelled \(a+b+1\) in \(D_i^u\). Now, denote with \(d_i = D_i^u = D_{i+s}^d\) and set \(\mathcal{D} = \{d_1, \dots, d_n\}, \mathcal{E} = \{e_1, \dots, e_m\}\).  

Naturally, if we cut \(\Sigma_n\) along \(\mathcal{D}\) we do not disconnect the surface; if \(m = n\) and even cutting \(\Sigma_n\) along \(\mathcal{E}\) does not disconnect the surface, then \((\Sigma_n, \mathcal{D}, \mathcal{E})\) is a Heegaard diagram of genus \(n\) of a 3-manifold, completely determined by the 6-tuple \((a,b,c,n,r,s)\): we call such a manifold \emph{Dunwoody manifold}.  

Clearly not all the 6-tuples \(\sigma=(a,b,c,n,s,r) \in \mathbb{Z}^6\) determine a Dunwoody manifold, that is,  are such that the set \(\mathcal{E}\) contains exactly \(n\) closed curves not disconnecting the surface $\Sigma_n$. We call \textit{admissible} such   6-tuples
and with  \(D(\sigma)\) (resp. \(H(\sigma)\)),  the   Dunwoody  manifold  (resp. open Heegaard diagram) associated to \(\sigma\in S\).

\begin{figure}[h!]
    \centering
    \includegraphics[width = .5\textwidth]{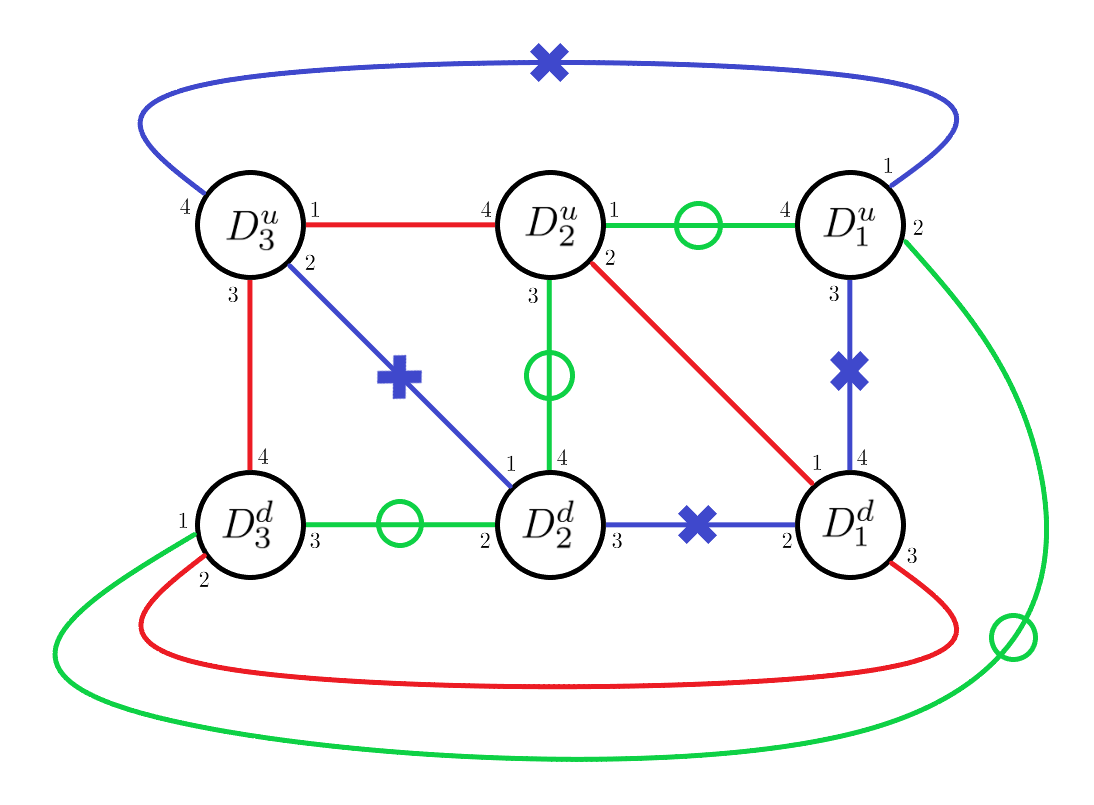   }
    \caption{The Dunwoody manifold  $M(1,1,1,3,2,1)$  homeomorphic to $S^1\times S^1\times S^1$: the three curves of the $\mathcal E$ system are depicted with different colors and different (or no) labels.}
    \label{fig:example_dunwoody}
\end{figure}

In Figure \ref{fig:example_dunwoody} we depict the case of $M(1,1,1,3,2,1)$ that is homeomorphic to $S^1\times S^1\times S^1$ (see \cite{cattabriga2004all}).

\subsection{The dictionary for  Dunwoody manifolds} 
\begin{figure}[h!]
    \centering
\begin{subfigure}{.45\textwidth}
    \centering
    \includegraphics[width = \textwidth]{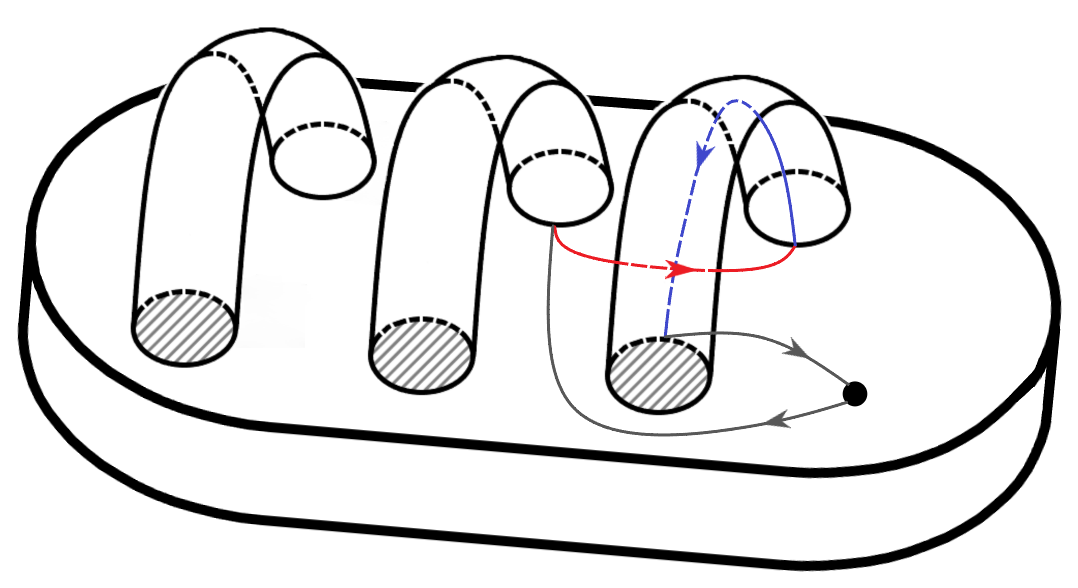}
\end{subfigure}
\begin{subfigure}{.45\textwidth}
    \centering
    \includegraphics[width = \textwidth]{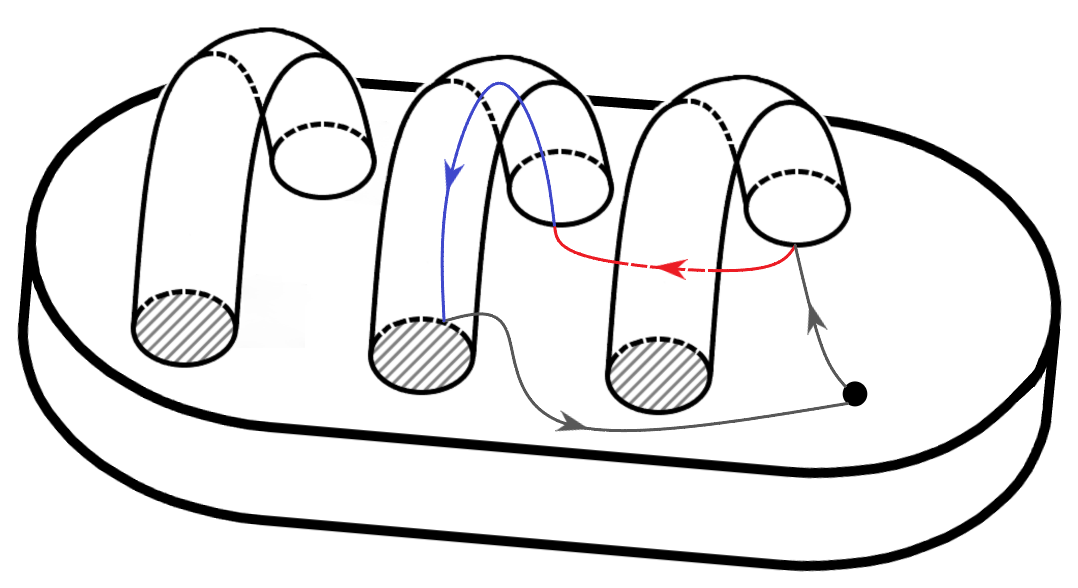}
\end{subfigure}

\begin{subfigure}{.45\textwidth}
    \centering
    \includegraphics[width = \textwidth]{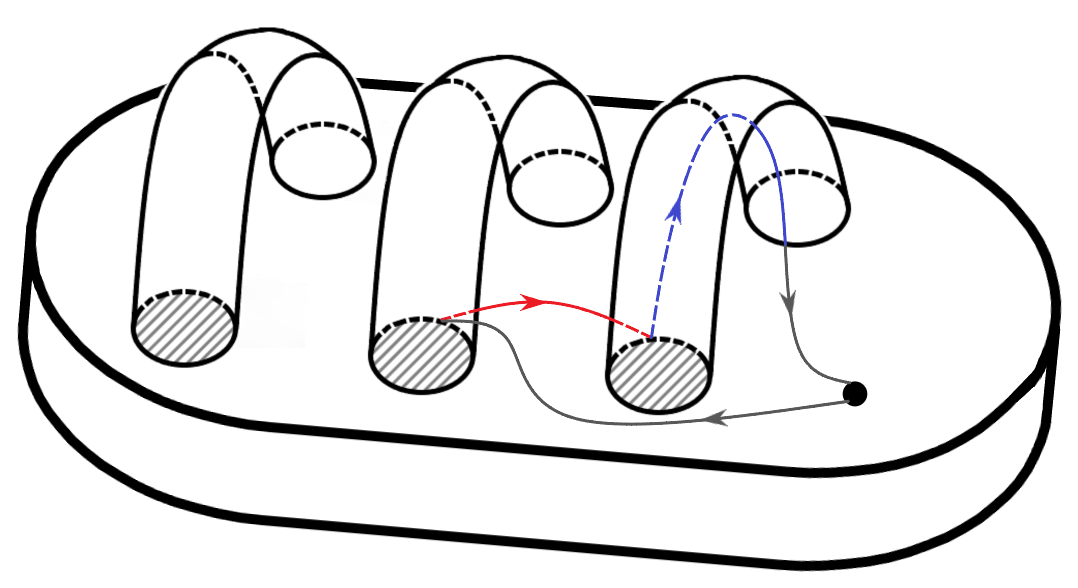}
\end{subfigure}
\begin{subfigure}{.45\textwidth}
    \centering
    \includegraphics[width = \textwidth]{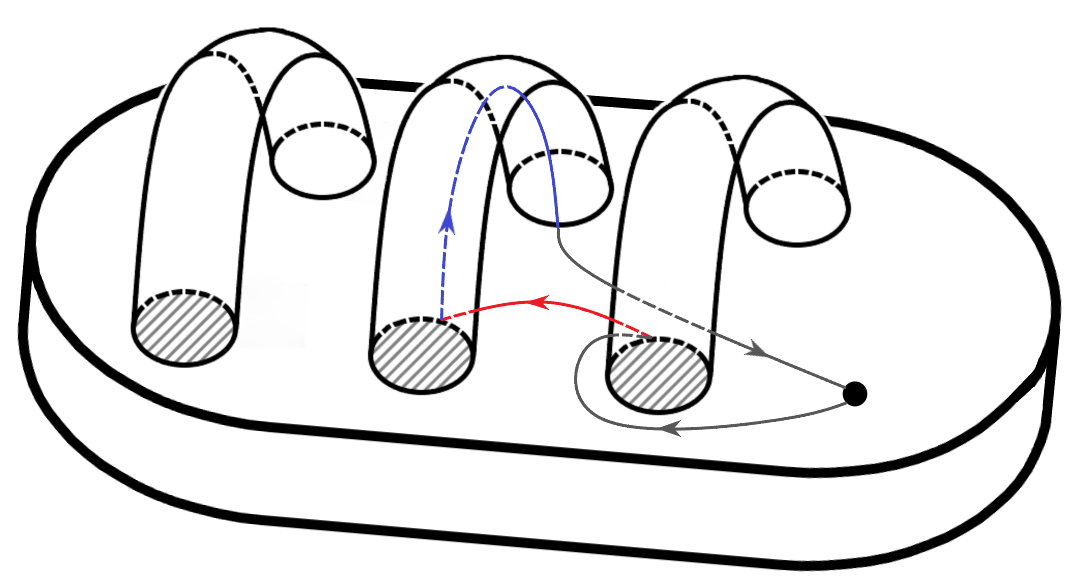}
\end{subfigure}
    \caption{The four possible horizontal arcs.}
    \label{fig:move_a}
\end{figure}

\begin{figure}[h!]
\centering
\begin{subfigure}{.45\textwidth}
    \centering
    \includegraphics[width = \textwidth]{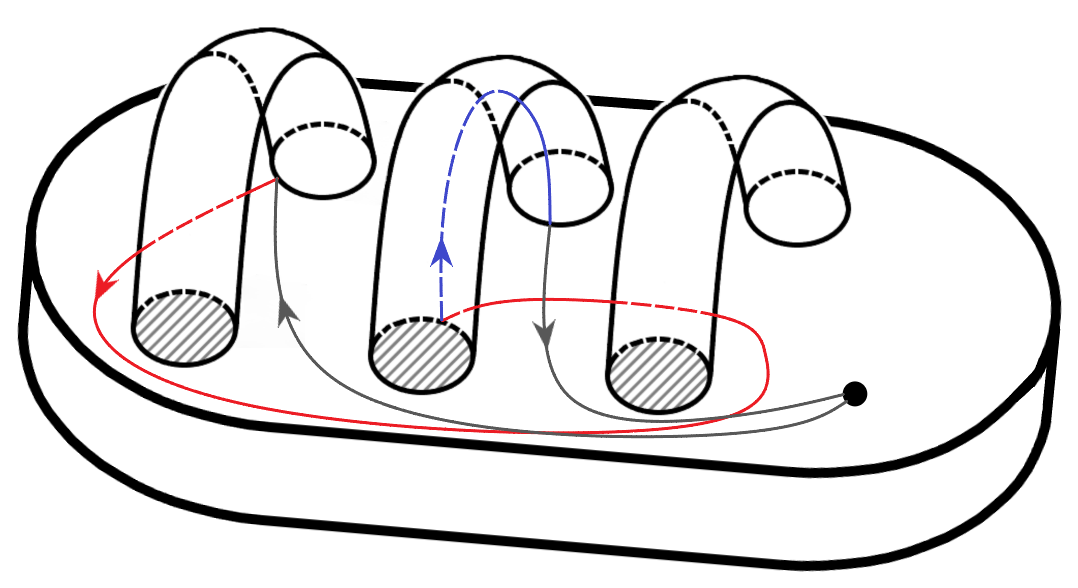}
\end{subfigure}
\begin{subfigure}{.45\textwidth}
    \centering
    \includegraphics[width = \textwidth]{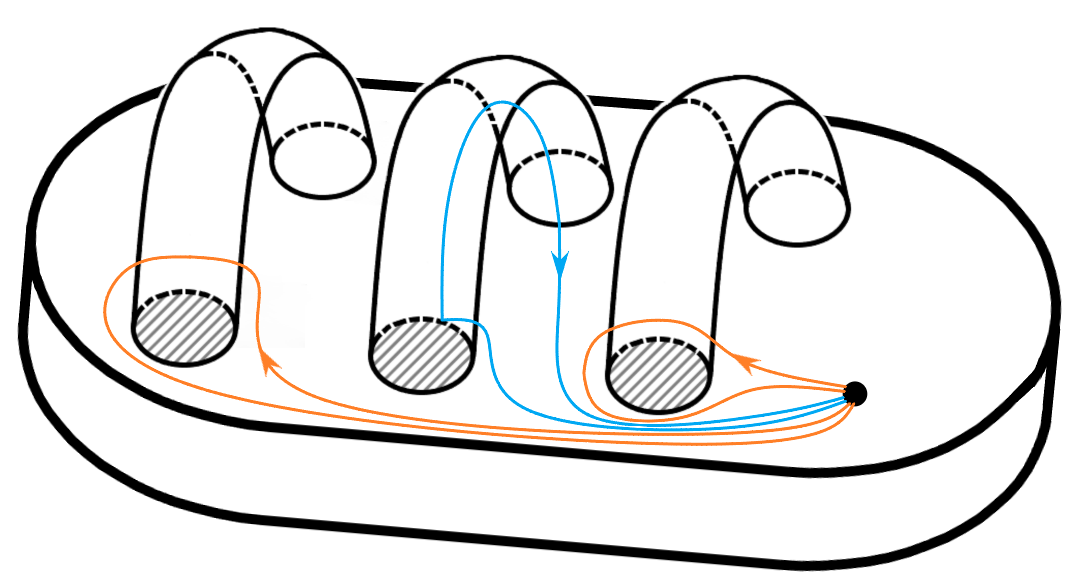}
\end{subfigure}
    \caption{An example of elementary curve corresponding to a vertical arc - in this case \(\downarrow \! C_3^3 = \beta_1 \beta_3 \alpha_2^{-1}\).}
    \label{fig:move_c}
\end{figure}

Given a Dunwoody manifold \(M=M(a,b,c,n,s,r)\), we use the open Heegaard diagram of Figure \ref{fig:dunwoody_label} for the construction of the dictionary. The elementary blue-red  pieces are  those corresponding to the four type of arcs of the graph, i.e., upper, lower, diagonal and vertical, each with two possible orientation. We use $P_1$ as base point (see Figure \ref{fig:generator_dunwoody}) and, in  the figures, we color in grey the fixed arcs connecting the elementary pieces with the base point.
\begin{itemize}
    \item Horizontal arc, upper or lower (see Figure \ref{fig:move_a}): we denote with $\overrightarrow{A^U_i}$ (resp. $\overrightarrow{A^L_i}$) an upper (resp. lower) arc connecting $D_i^u$ (resp. $D_i^d$) with $D_{i-1}^u$ (resp. $D_{i-i}^d$), and denote with $\overleftarrow{A^U_i}$ (resp. $\overleftarrow{A^L_i}$) an upper (resp. lower) arc connecting $D_i^u$ (resp. $D_i^d$) with $D_{i+1}^u$ (resp. $D_{i+i}^d$)
    \begin{align*}
        \overrightarrow{A^{U}_i} &= \beta_{i-1}^{-1} \alpha_{i-1} \\
        \overleftarrow{A^{U}_i} &= \beta_{i} \alpha_{i+1} \\
        \overrightarrow{A^{L}_i} &= \alpha_{i-s-1}^{-1} \\
        \overleftarrow{A^{L}_i} &= \alpha_{i-s+1}^{-1}.
    \end{align*}
    \item Vertical arc (see Figure \ref{fig:move_c}):  the description of the elementary piece corresponding to the vertical arc depends on the gluing parameter $s$ of 6-tuple determining $M$. We denote with  $C_i^s$ the elementary blue-red piece corresponding to a vertical arc connecting the disks $D_i^u$ and $D^d_i$ with parameter $s$ (so that $D^u_i$ is glued to $D^d_{i+s}$), while we use the arrow on the left to denote the orientation
        \begin{align*}
            \downarrow \! C_i^s &= w_{i,s} \alpha_{i+n-s}^{-1}\\
            \uparrow \! C_i^s &= w_{i,s}^{-1} \alpha_{i}
        \end{align*}
        where 
         \[w_{i,s} = \prod_{j = 0}^{k} \beta_{i+j}, \quad 1\leq k\leq n \textup{ and } k\equiv n-s-1 \textup{ mod } n.\]
    \item Diagonal arc: from the construction  of the Dunwoody manifold it is quite straightforward to see that the elementary curve corresponding to a diagonal arc with gluing parameter $s$ is equal to the one corresponding to a vertical arc with gluing parameter $s+1$, that is
        \[\downarrow \! B_i^s = \downarrow \! C_i^{s+1}\]
        \[\uparrow \! B_{i}^s = \uparrow \! C_{i+1}^{s+1},\]
        where $\downarrow \!B_i^s$ denotes the elementary piece corresponding to a diagonal arc connecting  the disks $D_i^u$ and $D^d_{i-1}$ with parameter $s$, and $\uparrow \! B_i^s$ denotes the elementary piece corresponding to a diagonal arc connecting  the disks $D_i^d$ and $D^u_{i+1}$ with parameter $s$.
\end{itemize}

\begin{example}
In the case of the Dunwoody manifold  $M(1,1,1,3,2,1)$, homeomorphic to $S^1\times S^1\times S^1$ (see \cite{cattabriga2004all}) and depicted in Figure \ref{fig:example_dunwoody}, if we denote with $e_1$ the blue curve (label $\times$), $e_2$ the green one (label $\circ$) and $e_3$ the red one (no label) we have, starting from the vertex labeled 3:
\[\begin{array}{ll}
\overline{e_1} &= \beta_1 \beta_2 \alpha_3^{-1} \beta_3 \alpha_1 \beta_3^{-1} \alpha_3 \alpha_1^{-1}\\
\overline{e_2} &= \beta_2 \beta_3 \alpha_1^{-1} \beta_1 \alpha_2 \beta_1^{-1} \alpha_1 \alpha_2^{-1}\\ 
\overline{e_3} &= \beta_3 \beta_1 \alpha_2^{-1} \beta_2 \alpha_3 \beta_2^{-1} \alpha_2 \alpha_3^{-1}.
\end{array}\]
\end{example}

\subsection{Fibonacci and Sieradski manifolds} 
Recall that the Minkus manifold  $M_n(2a+1,2r)$ is the $n$-fold cyclic covering branched on the 2-bridge knot $b(2a+1,2r)$, with $(2a+1, 2r) = 1$. In   \cite{grasselli2001genus} it is proved that $M_n(2a+1,2r) = M(a,0,1,n,r,\bar{s})$, where $\bar{s}$ can be computed on $\Gamma(a,0,1,n,r,0)$ as follows: consider the vertex labelled \(v = a+b+1\) of \(D^u_1\), and orient downwards the arc connecting this vertex with \(D^d_1\) and call \(e_1\) the cycle of \(\mathcal E\) passing from $v$. Now, follow along this orientation the arcs of \(\Gamma\) belonging to \(e_1\), and count the arcs running from  $D^u_{i-1}$ to $D^u_{i}$ or $D^d_{i-1}$ to $D^d_{i}$ and those running from  $D^u_{i}$ to $D^u_{i-1}$ or $D^d_{i}$ to $D^d_{i-1}$. The difference between the first type of arcs and the second one is \(\bar{s}\). 

\begin{figure}[h!]
    \centering
    \includegraphics[width = .5\textwidth]{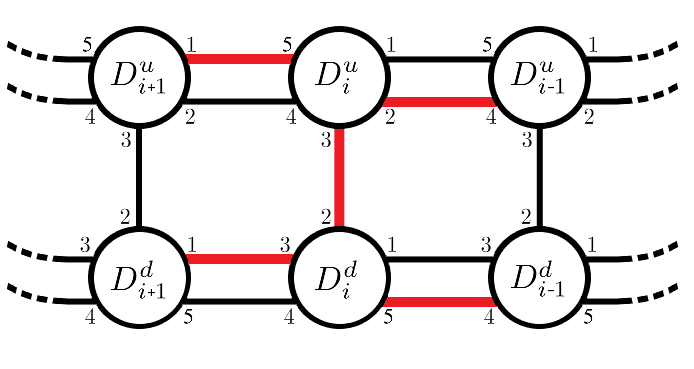}
    \caption{The general case of a Fibonacci manifold, described by the 6-tuple $M(2,0,1,n,1,0)$}
    \label{fig:fibonacci_general}
\end{figure}

\begin{figure}[h!]
    \centering
    \includegraphics[width = .8\textwidth]{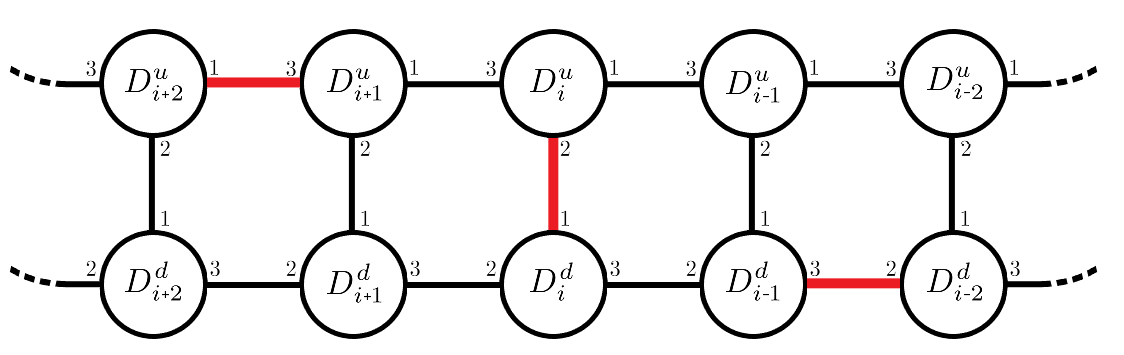}
    \caption{The general case of a Sieradski manifold, described by the 6-tuple $M(1,0,1,n,1,-2), n>2$}
    \label{fig:sieradski_general}
\end{figure}

For $r=1$ and $a=2$, the knot $b(5,2)$ is the figure-eight knot, so $M(2,0,1,n,1,\bar{s})$ is a Fibonacci manifold, while when $r=a=1$, the knot $b(3,2)$ is the trefoil knot, so $M(1,0,1,n,1,\bar{s})$ is a Sieradski manifold.
For both families of manifolds it is possible to compute explicitly the words related to the curves: 
\begin{itemize}
    \item Fibonacci case $M(2,0,1,n,1,\bar{s})$: we have $\bar{s} = 0$ fal all $n \in \mathbb{N}$; so the curve $e_i$ starting in the vertex labeled 3 of \(D_i^u\) is the red curve (and thicker curve) in Figure \ref{fig:fibonacci_general}, with related word: 
    \[\overline{e_i} = \alpha_{i}^{-1} \beta_{i-1}^{-1} \alpha_{i-1} \alpha_{i}^{-1} \beta_{i} \alpha_{i+1} \alpha_i^{-1}.\]
    \item Sieradski case $M(1,0,1,n,1,\bar{s})$: we have $\bar{s} = -2$ for $n>2$; so the curve $e_i$ starting in the vertex labeled 2 of \(D_i^u\) is the red curve (and thicker curve) in Figure \ref{fig:sieradski_general}, with related word: 
    \[\overline{e_i} = \beta_{i} \beta_{i+1} \alpha_{i+2}^{-1} \beta_{i+1}^{-1} \alpha_{i+1} \alpha_{i}^{-1}.\]
\end{itemize}

\section{Periodic Takahashi manifolds}\label{Takahashi}
\begin{figure}[h!]
 \centering
 \includegraphics*[width = .5\textwidth]{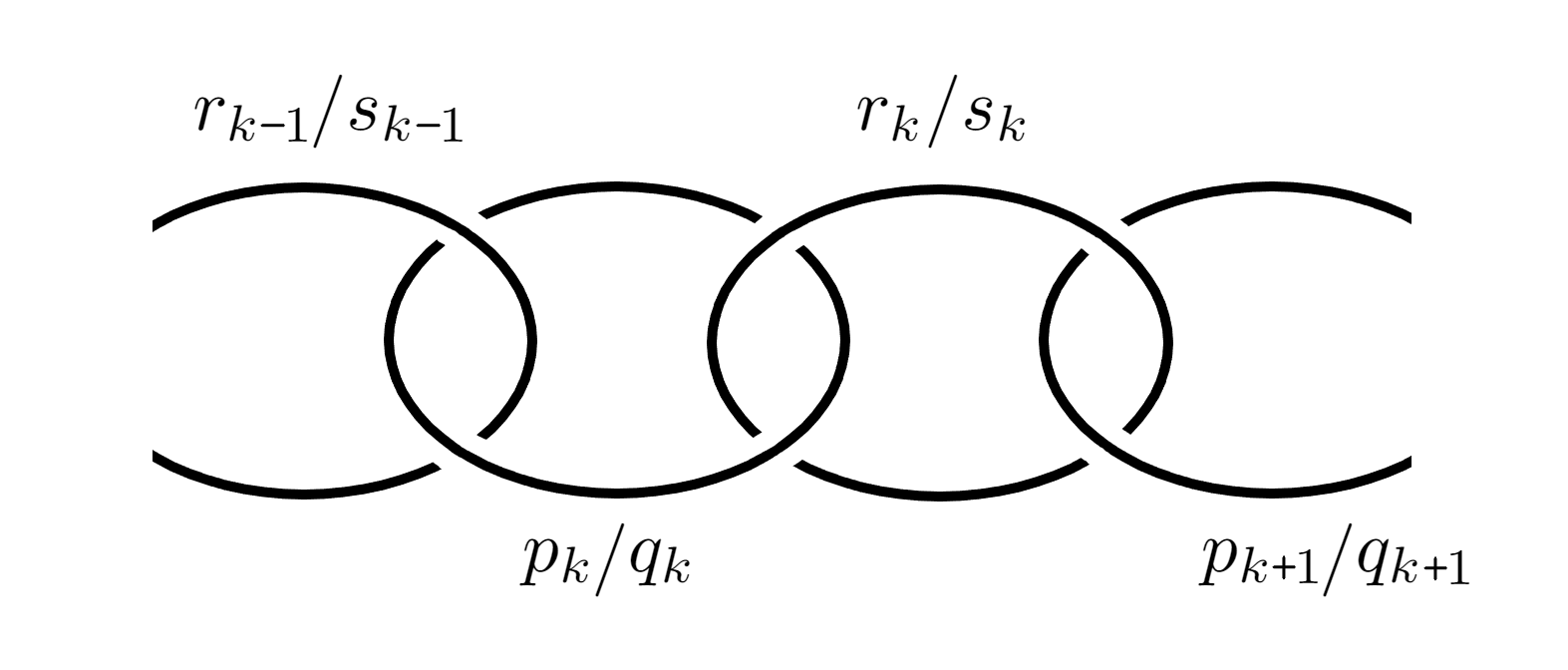}
 \caption{Presentation via Dehn surgery for the Takahashi manifold   $T_{n} (p_1 / q_1,\ldots, p_n/q_n; r_1/s_1,\ldots, r_n/s_n)$.}
 \label{fig:takahashi_dehn_basic}
\end{figure}

Let  ${\cal L}_{n} \subseteq  S^3$ be the link with $2n$ components represented in Figure \ref{fig:takahashi_dehn_basic}: it is a closed chain of $2n$ unknotted components having a cyclic symmetry of order $n$ which permutes the unknotted components.  

The \textit{Takahashi manifold}  $T_{n} (p_1 / q_1,\ldots, p_n/q_n; r_1/s_1,\ldots, r_n/s_n)$ is the manifold obtained by Dehn surgery on $S^3$, along the link ${\cal L}_{n}$, with surgery coefficients alternatively equal to $p_{k} / q_{k}$  and $r_{k} / s_{k}$, with $1 \le k \le n$. Without loss of generality, we can always suppose $\gcd(p_{k},q_{k})=1$, $\gcd(r_{k},s_{k})=1$ and $p_{k},r_{k}\ge 0$. If $p_1/q_1=p_2/q_2=\cdots =p_m/q_m=p/q$ and $r_1/s_1=r_2/q_2=\cdots =r_m/s_m=r/s$, i.e., the surgery coefficients have the same cyclic symmetry of order $n$ of the link ${\cal L}_{n}$, the Takahashi manifold is called \textit{periodic} and it is denoted simply by $T_n(p/q,r/s)$.

\begin{figure}[h!]
    \centering
\begin{subfigure}{.45\textwidth}
    \centering
    \includegraphics[width = \textwidth]{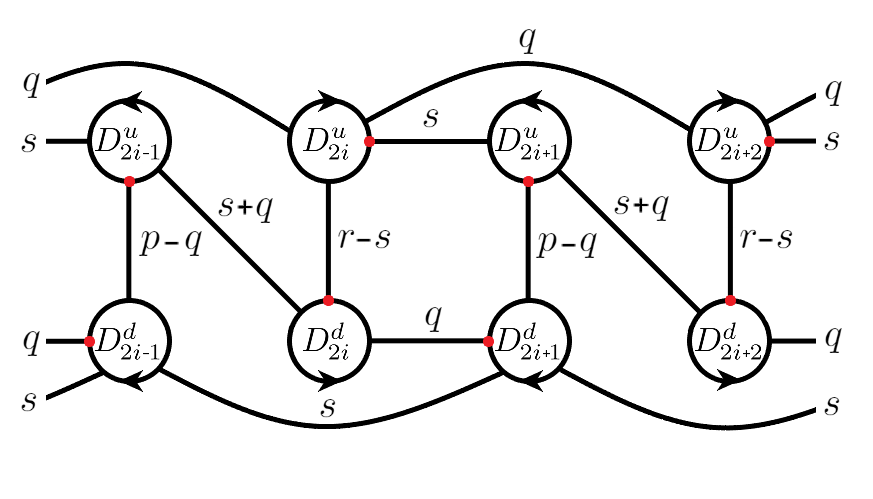}
    \caption{\(p>q, r>s\)}
\end{subfigure}
\begin{subfigure}{.45\textwidth}
    \centering
    \includegraphics[width = \textwidth]{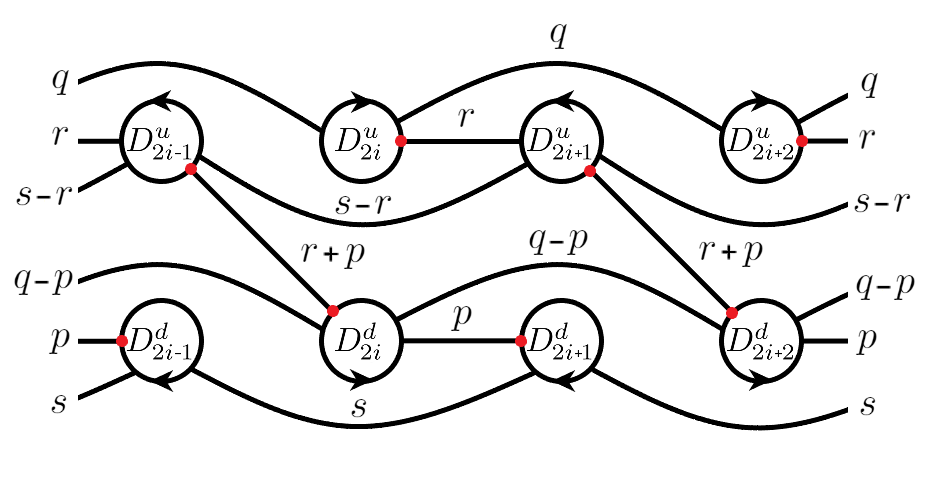}
    \caption{\(p<q, r<s\)}
\end{subfigure}

\begin{subfigure}{.45\textwidth}
    \centering
    \includegraphics[width = \textwidth]{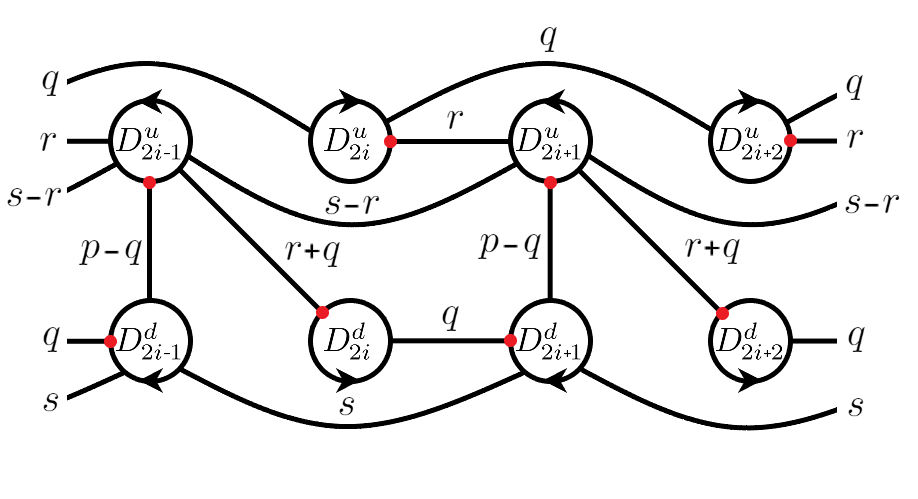}
    \caption{\(p>q, r<s\)}
\end{subfigure}
\begin{subfigure}{.45\textwidth}
    \centering
    \includegraphics[width = \textwidth]{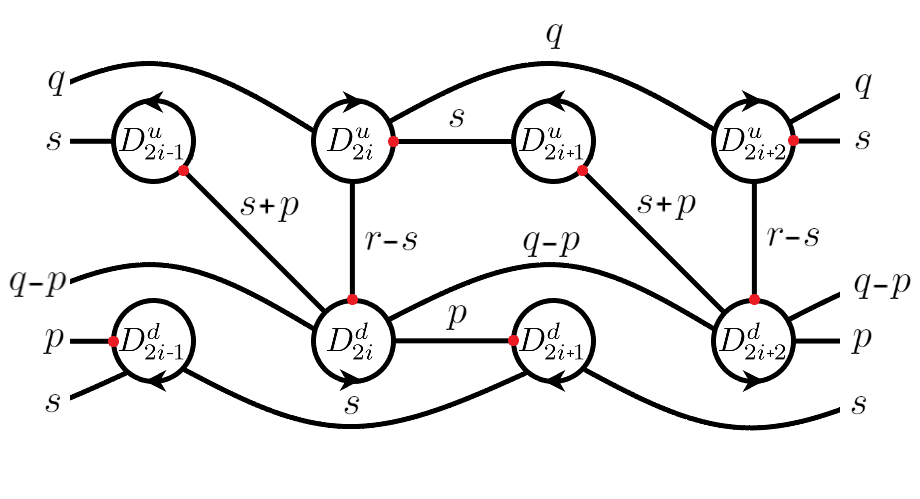}
    \caption{\(p<q, r>s\)}
\end{subfigure}
 \caption{Four open Heegaard diagrams describing, depending on the parameters \(p,q,r,s\), the periodic  Takahashi manifold. An arc with label \(k\) identifies \(k\) parallel arcs.}
 \label{fig:takahashi_heeg_2}
\end{figure}

In \cite[Proposition 7]{cristofori2007strongly} the authors describe an open Heegaard diagram for a periodic  Takahashi manifold   $T_n(p/q,r/s)$ in the case $p/q,r/s\geq 0$ and $p/q,r/s\ne 1$, that, up to isotopy, is  the one depicted in  Figure \ref{fig:takahashi_heeg_2}, where the couples of corresponding circles are \(D_{i}^u, D_{i}^d\) for \(i = 1, \dots, 2n\) glued according to the orientation and so that the fat red points are identified; without loss of generality we always assume $p,q,r,s\geq 0$.  In Figure \ref{fig:example_takahashi} we depict the case of $T_2(1/2, 2/3)$.


 

\subsection{The dictionary for periodic Takahashi manifolds}

\begin{figure}[h!]
    \centering
    \includegraphics[width = .6\textwidth]{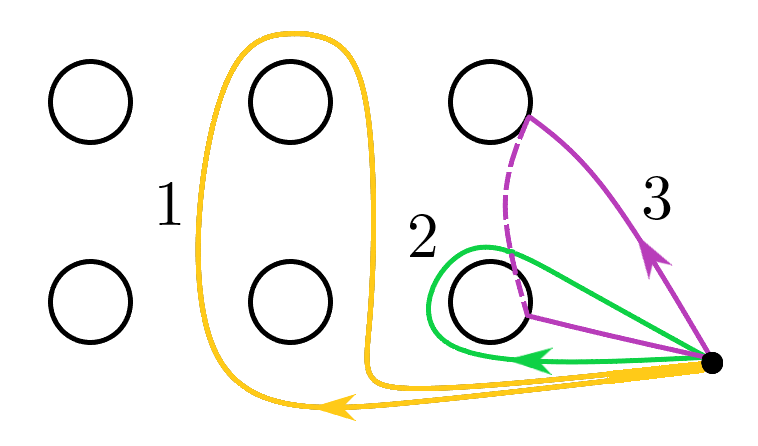}
    \caption{The case of an arc of type \(\overrightarrow{F}\).}
    \label{fig:taka_resol_c}
\end{figure}

For the construction of the dictionary in the setting of Takahashi manifolds, we refer to the Heegaard diagram of Figure \ref{fig:takahashi_heeg_2}: we have eight  type  of elementary arcs and we denote them according to Figure \ref{fig:taka_arcs}. Following Figure \ref{fig:taka_moves} down below, we can describe the words associated to the  elementary red-blue arcs. As in the Dunwoody case, we use $P_1$ as the base point and, in  the figures, we color in grey the  arcs connecting the elementary pieces with the base point. 
\begin{itemize}
    \item Arcs \(A^U\) and \(A^L\): they are as the horizontal arcs in the Dunwoody case  so the result is the same (see Figure \ref{fig:move_a}) 
        \begin{align*}
            \overrightarrow{A^{U}_{2i}} &= \beta_{2i+1}^{-1} \alpha_{2i+1} \\
            \overleftarrow{A^{U}_{2i}} &= \beta_{2i+1} \alpha_{2i} \\
            \overrightarrow{A^{L}_{2i}} &= \alpha_{2i+1}^{-1} \\
            \overleftarrow{A^{L}_{2i}} &= \alpha_{2i}^{-1}.
        \end{align*}
    \item Arcs \(B\) and \(C\): these coincide with vertical or diagonal arcs of the Dunwoody case with \(s = 0\), so we have
        \begin{align*}
            \overrightarrow{B_{2i-1}} &= \alpha_{2i}^{-1} \\
            \overleftarrow{B_{2i-1}} &= \alpha_{2i-1} \\ 
            \downarrow \! C_{2i-1} &= \alpha_{2i-1}^{-1} \\
            \uparrow \! C_{2i-1} &= \alpha_{2i-1}. 
        \end{align*}
    \item Arc \(F\): we obtain the curve depicted in Figure \ref{fig:taka_moves}. In order to find the corresponding word we need to decompose the original loop into three elements, see Figure \ref{fig:taka_resol_c}. Following these three elementary pieces we obtain
        \begin{align*}
            \overrightarrow{F_{2i}} &= [\alpha_{2i+1}, \beta_{2i+1}] \beta_{2i+2}^{-1} \alpha_{2i+2} = \alpha_{2i+1}^{-1} \beta_{2i+1}^{-1} \alpha_{2i+1} \beta_{2i+1} \beta_{2i+2}^{-1} \alpha_{2i+2} \\
            \overleftarrow{F_{2i}} &= \beta_{2i} [\beta_{2i-1}, \alpha_{2i-1}] \alpha_{2i-2} = \beta_{2i} \beta_{2i-1}^{-1} \alpha_{2i-1}^{-1} \beta_{2i-1} \alpha_{2i-1} \alpha_{2i-2}
        \end{align*}
    where \([\alpha,\beta] = \alpha^{-1} \beta^{-1} \alpha \beta\). 
    \item Arcs \(G\), \(X\) and \(Y\): for all these kind of arcs the procedure is similar to the previuos one, obtaining the following sequence of moves
        \begin{align*}
            \overrightarrow{G_{2i-1}} &= \beta_{2i}^{-1} \beta_{2i+1}^{-1} \alpha_{2i+1} \\
            \overleftarrow{G_{2i-1}} &= \beta_{2i-1} \beta_{2i-2} \alpha_{2i-3} \\
            \overrightarrow{X_{2i-1}} &= \beta_{2i+1}^{-1} \alpha_{2i+2}^{-1} \\
            \overleftarrow{X_{2i-1}} &= \beta_{2i-1}^{-1} \alpha_{2i-2}^{-1} \\
            \overrightarrow{Y_{2i-1}} &= \alpha_{2i+1}^{-1} \\
            \overleftarrow{Y_{2i-1}} &= \alpha_{2i-3}^{-1}.
        \end{align*}
\end{itemize}

\begin{figure}[h!]
    \centering
    \includegraphics[width = .5\textwidth]{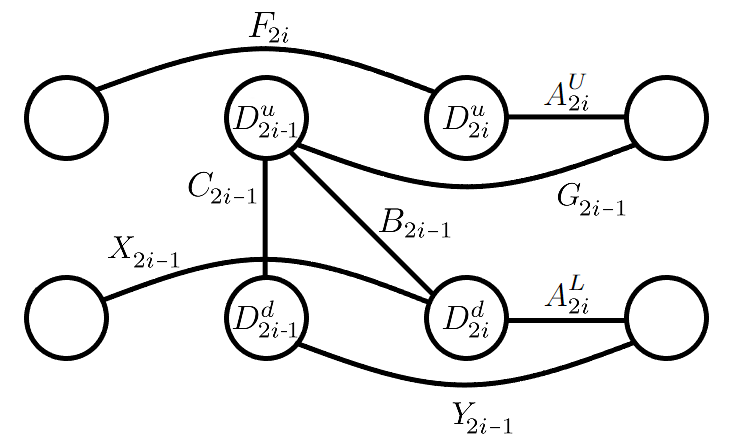}
    \caption{All the possible elementary arcs in the periodic Takahashi case.} 
    \label{fig:taka_arcs}
\end{figure}

\begin{figure}[h!]
    \centering
\begin{subfigure}{.45\textwidth}
    \centering
    \includegraphics[width = \textwidth]{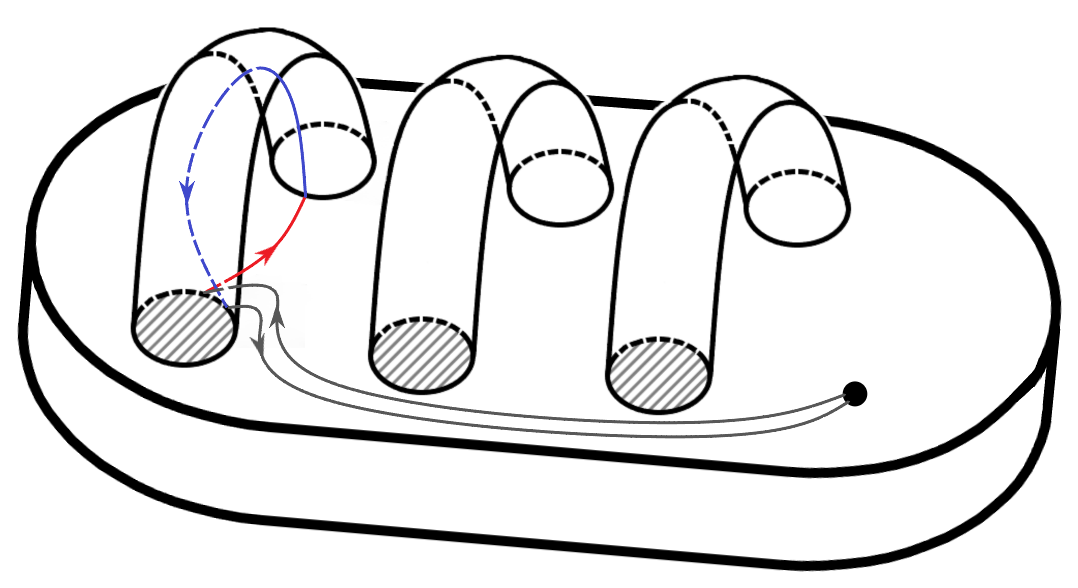}
\end{subfigure}
\begin{subfigure}{.45\textwidth}
    \centering
    \includegraphics[width = \textwidth]{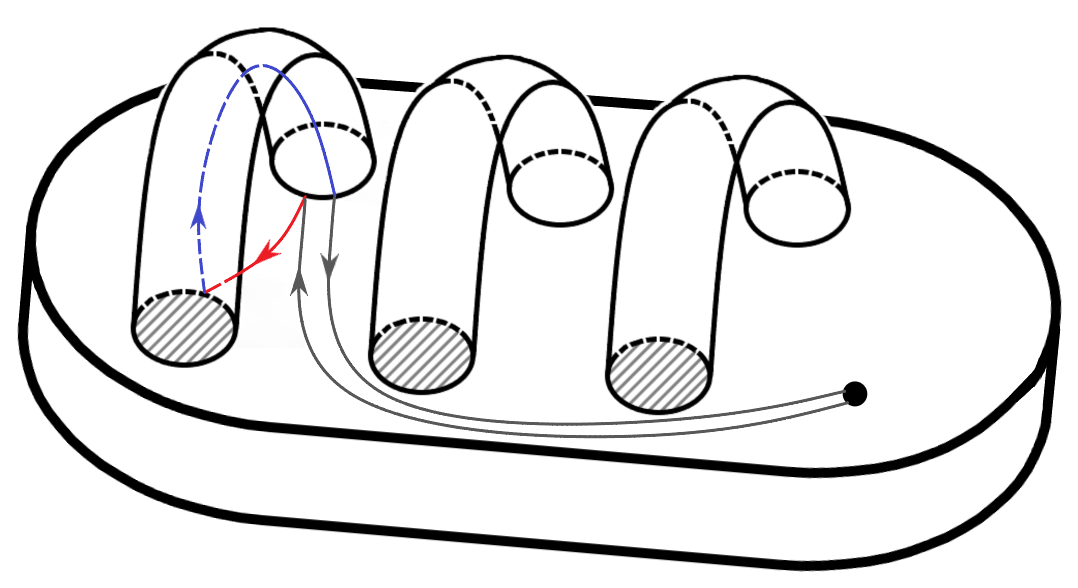}
\end{subfigure}

\begin{subfigure}{.45\textwidth}
    \centering
    \includegraphics[width = \textwidth]{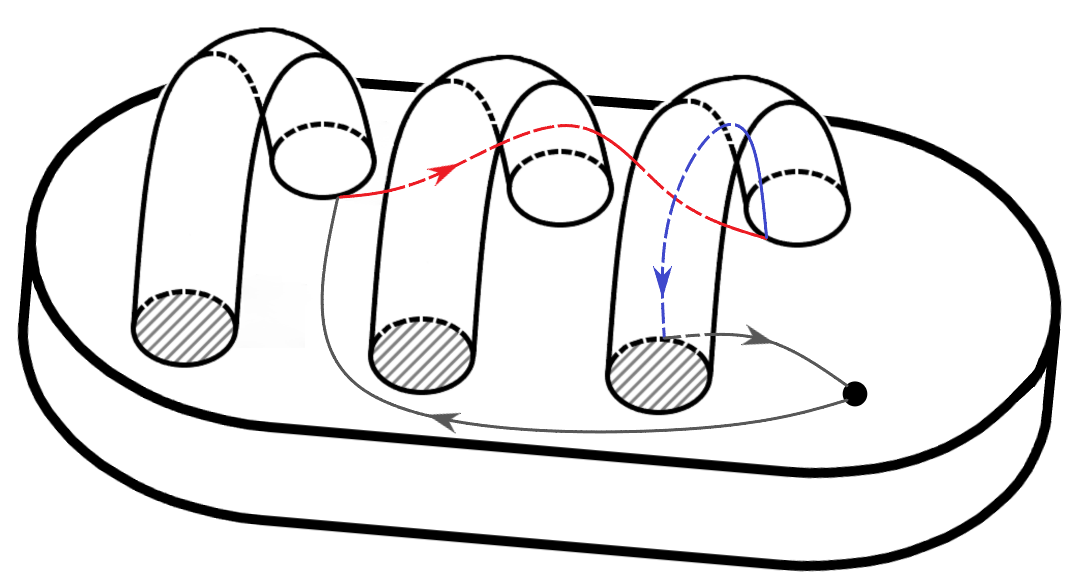}
\end{subfigure}
\begin{subfigure}{.45\textwidth}
    \centering
    \includegraphics[width = \textwidth]{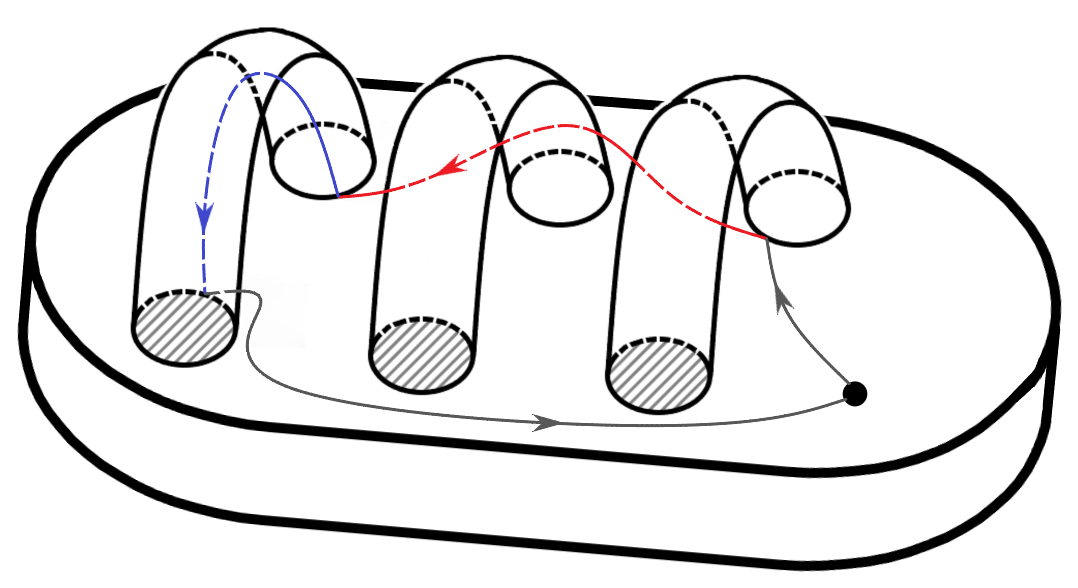}
\end{subfigure}

\begin{subfigure}{.45\textwidth}
    \centering
    \includegraphics[width = \textwidth]{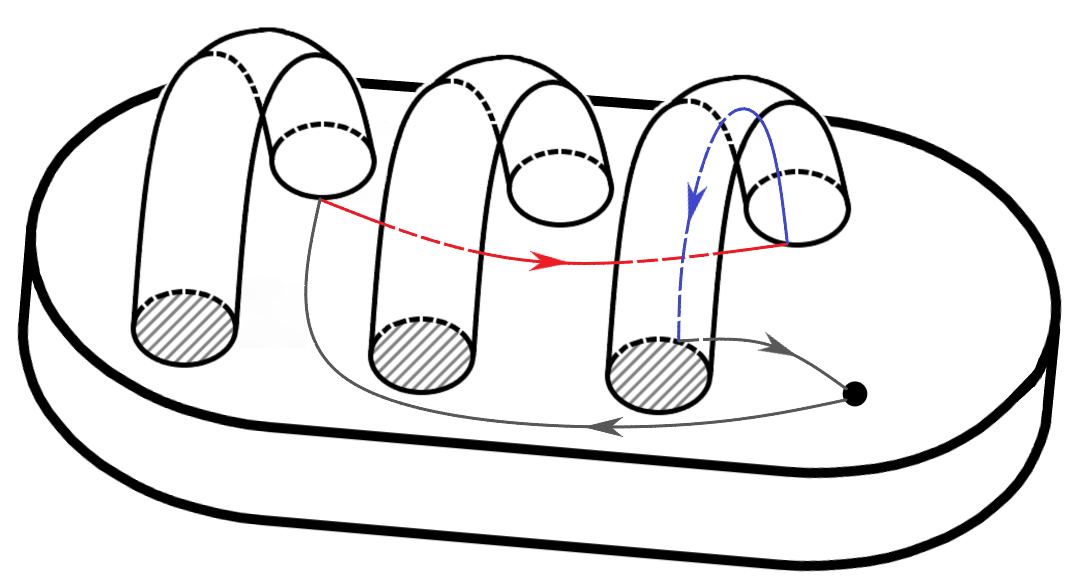}
\end{subfigure}
\begin{subfigure}{.45\textwidth}
    \centering
    \includegraphics[width = \textwidth]{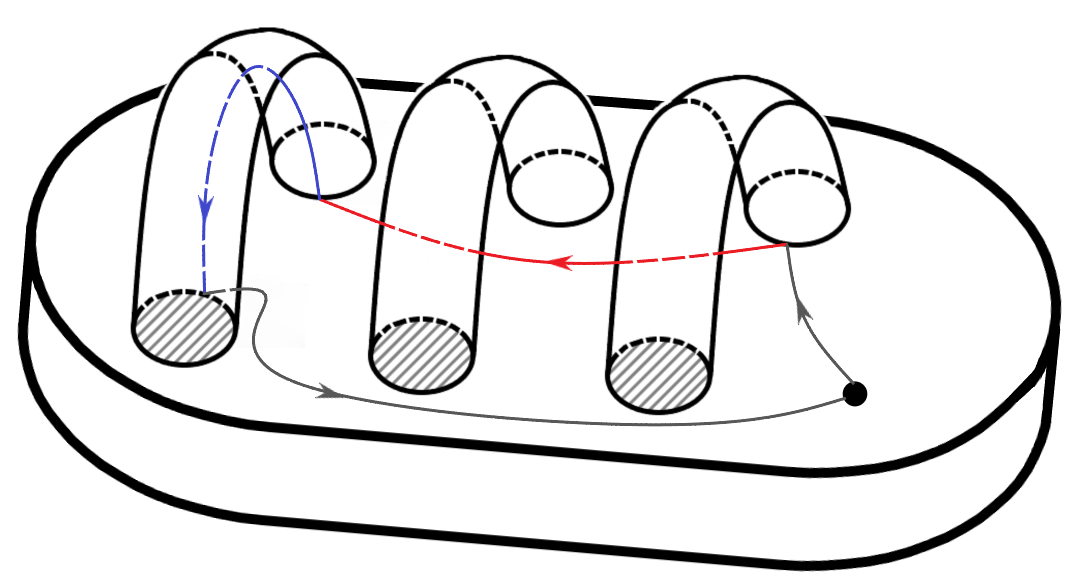}
\end{subfigure}

\begin{subfigure}{.45\textwidth}
    \centering
    \includegraphics[width = \textwidth]{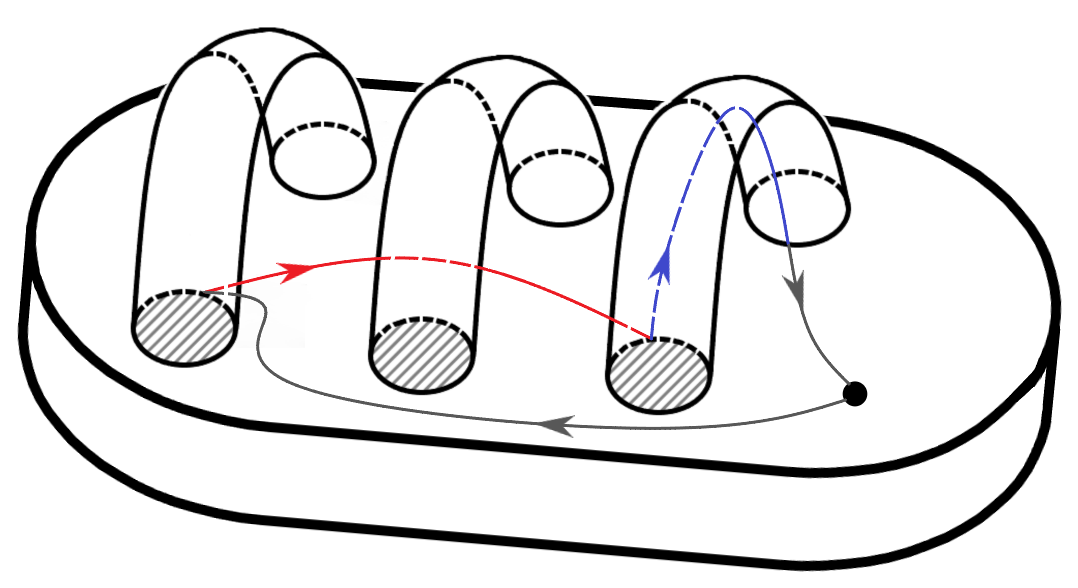}
\end{subfigure}
\begin{subfigure}{.45\textwidth}
    \centering
    \includegraphics[width = \textwidth]{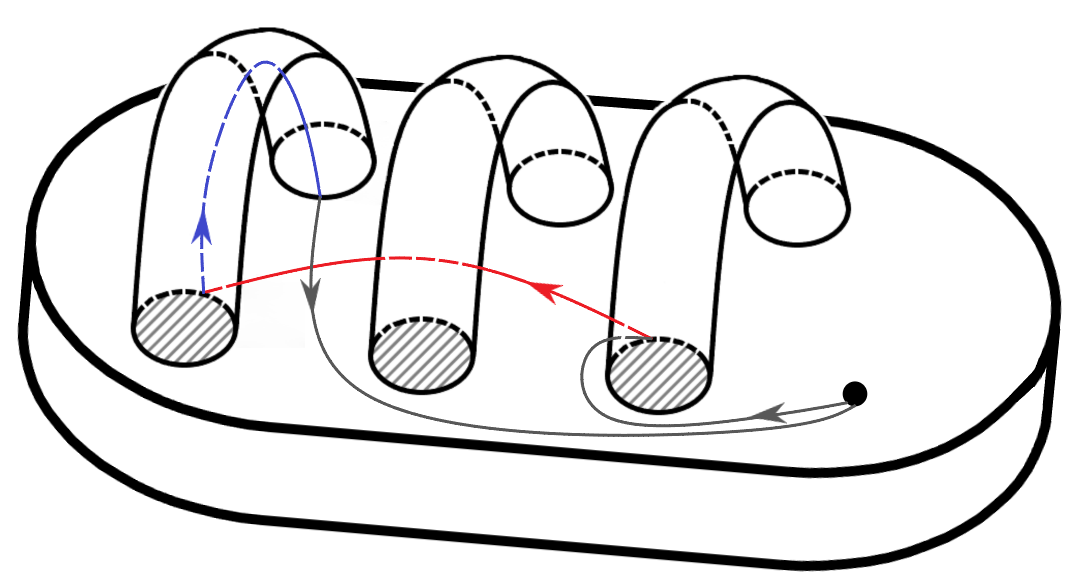}
\end{subfigure}

\begin{subfigure}{.45\textwidth}
    \centering
    \includegraphics[width = \textwidth]{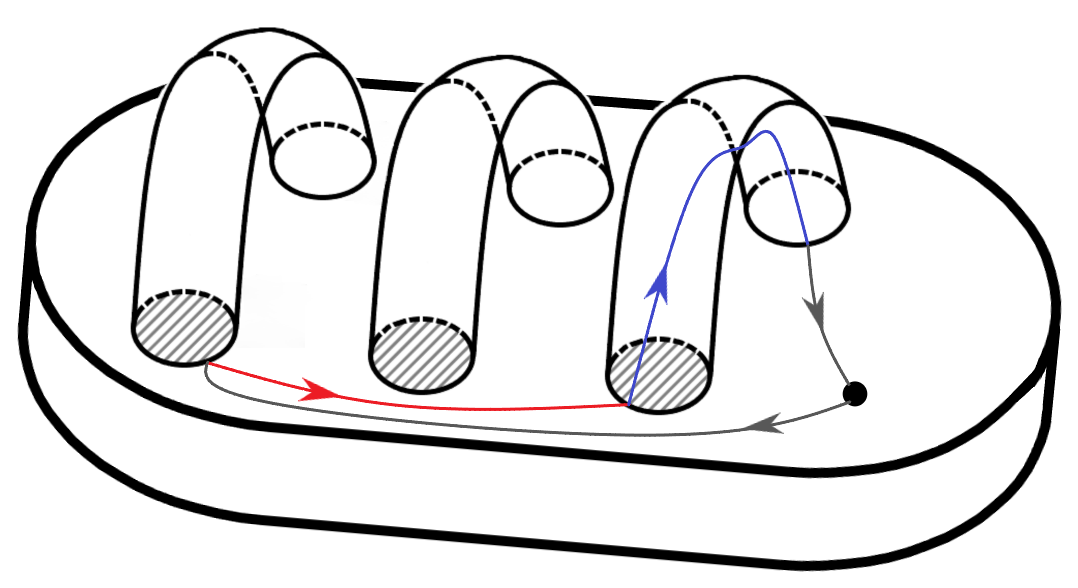}
\end{subfigure}
\begin{subfigure}{.45\textwidth}
    \centering
    \includegraphics[width = \textwidth]{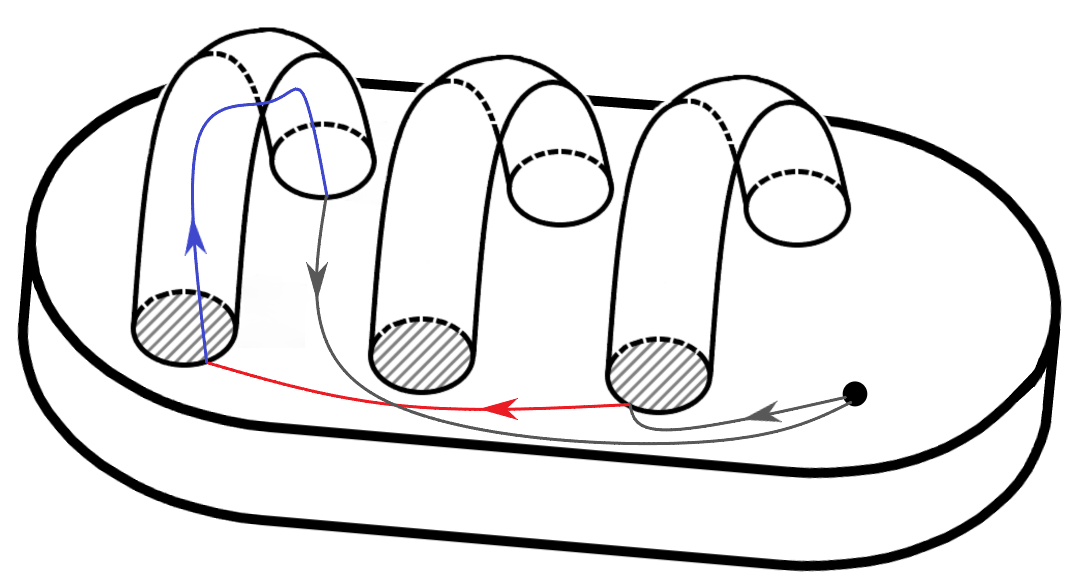}
\end{subfigure}

\begin{subfigure}{.45\textwidth}
    \centering
    \includegraphics[width = \textwidth]{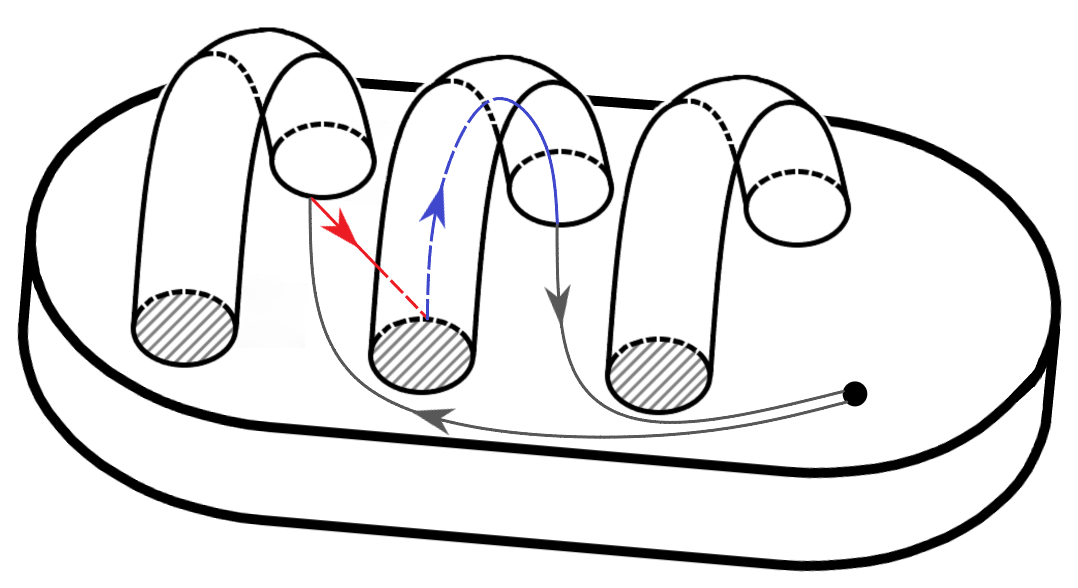}
\end{subfigure}
\begin{subfigure}{.45\textwidth}
    \centering
    \includegraphics[width = \textwidth]{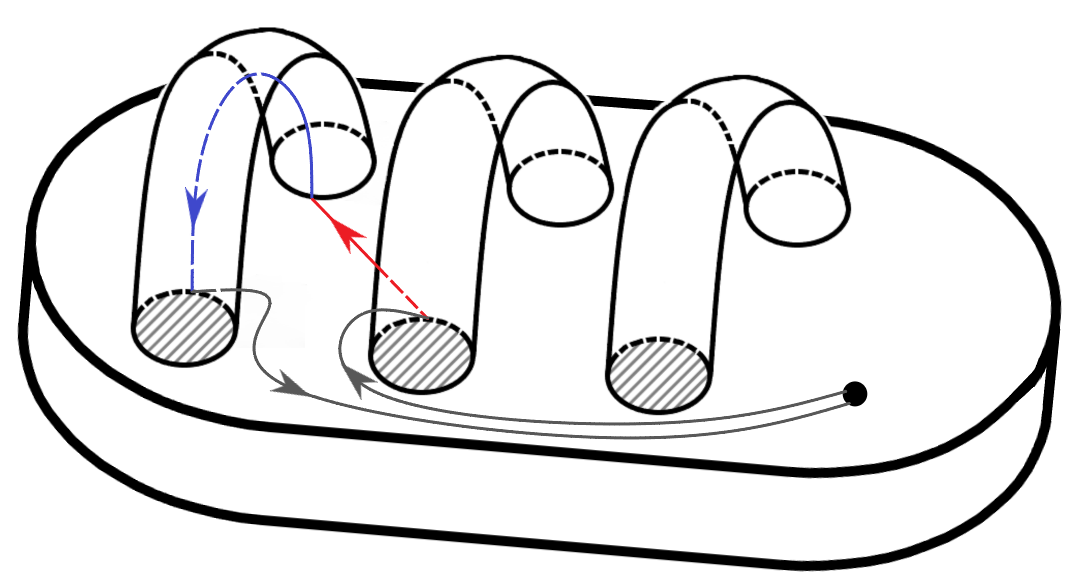}
\end{subfigure}

\caption{From top to bottom, the curves corresponding to  \(B\)-arcs, \(C\)-arcs, \(D\)-arcs, \(E\)-arcs, \(F\)-arcs and \(G\)-arcs.}
\label{fig:taka_moves}
\end{figure}

\begin{figure}[h!]
    \centering
    \includegraphics[width = .7\textwidth]{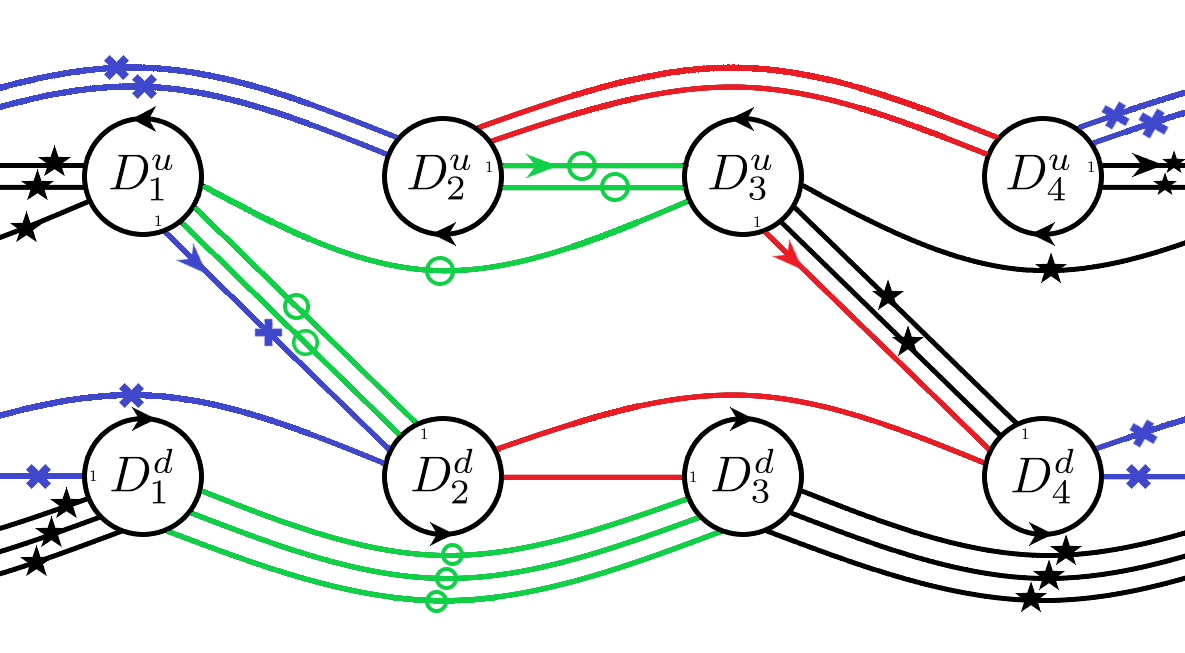}
    \caption{The Takahashi manifold \(T_2(1/2, 2/3)\).}
    \label{fig:example_takahashi}
\end{figure}

\begin{example}
In the case of the Takahashi manifold  $T_2(1/2, 2/3)$  depicted in Figure \ref{fig:example_takahashi}, if we denote with $e_1$ the blue curve (label $\times$), $e_2$ the red one (no label), $e_3$ the green one (label $\circ$) and $e_4$ the black one (label $\star$), following the arrow we have
\[\begin{array}{lll}
\overline{e_1} &= \overrightarrow{B_1} \overleftarrow{F_2} \overrightarrow{X_3} \overleftarrow{F_2} \overrightarrow{A^L_4} &= \alpha_2^{-1} \beta_2 \beta_1^{-1} \alpha_1^{-1} \beta_1 \alpha_1 \alpha_4 \beta_1^{-1} \alpha_2^{-1} \beta_2 \beta_1^{-1} \alpha_1^{-1} \beta_1 \alpha_1 \alpha_4 \alpha_1^{-1}\\
\overline{e_2} &= \overrightarrow{B_3} \overleftarrow{F_4} \overrightarrow{X_1} \overleftarrow{F_4} \overrightarrow{A^L_2} &= \alpha_4^{-1} \beta_4 \beta_3^{-1} \alpha_3^{-1} \beta_3 \alpha_3 \alpha_2 \beta_3^{-1} \alpha_4^{-1} \beta_4 \beta_3^{-1} \alpha_3^{-1} \beta_3 \alpha_3 \alpha_2 \alpha_3^{-1}\\ 
\overline{e_3} &= \overrightarrow{A^U_2} \overleftarrow{Y_3} \overrightarrow{G_1} \overleftarrow{Y_3} \overrightarrow{B_1} \overrightarrow{A^U_2} \overleftarrow{Y_3} \overrightarrow{B_1} &= \beta_3^{-1} \alpha_3 \alpha_1^{-1} \beta_2^{-1} \beta_3^{-1} \alpha_3 \alpha_1^{-1} \alpha_2^{-1} \beta_3^{-1} \alpha_3 \alpha_1^{-1} \alpha_2^{-1}\\
\overline{e_4} &= \overrightarrow{A^U_4} \overleftarrow{Y_1} \overrightarrow{G_3} \overleftarrow{Y_1} \overrightarrow{B_3} \overrightarrow{A^U_4} \overleftarrow{Y_1} \overrightarrow{B_3} &= \beta_1^{-1} \alpha_1 \alpha_3^{-1} \beta_4^{-1} \beta_1^{-1} \alpha_1 \alpha_3^{-1} \alpha_4^{-1} \beta_1^{-1} \alpha_1 \alpha_3^{-1} \alpha_4^{-1}.
\end{array}\]
\end{example}


\bibliographystyle{acm}
\bibliography{main}

\end{document}